\documentclass [12pt,a4paper,reqno]{amsart}
 \input amssymb.sty
\newcommand{\ef}{\end{equation}}
\chardef\bslash=`\\ 

\newtheorem{thm}{Theorem} [section]
\newtheorem*{thm*}{Theorem}
\newtheorem{cor}[thm]{Corollary}
\newtheorem{lem}[thm]{Lemma}

\newtheorem{defn}[thm]{Definition}

\newtheorem{example}[thm]{Example}
\newtheorem{acknowledgment*}[thm] {Acknowledgment}

 \theoremstyle{remark}

 \newtheorem{remark}[thm]{Remark}

\newcommand{\thmref}[1]{Theorem~\ref{#1}}

\newcommand{\lemref}[1]{Lemma~\ref{#1}}

 \renewcommand{\sectionmark}[1]{}

\newcommand{\iy}{\infty}

\newcommand{\doe}{\overset{\text{def}}{=}}
 
\newcommand{\loc} {\operatorname{loc}}

 \date{}

\begin{document}

\title[A criteria for correct solvability in $L_p(\mathbb R)$]
{A criteria for correct solvability in $L_p(\mathbb R)$\\ of a
general Sturm-Liouville  equation}
\author[N.A. Chernyavskaya]{N.A. Chernyavskaya}
\address{Department of Mathematics and Computer Science, Ben-Gurion
University of the Negev, P.O.B. 653, Beer-Sheva, 84105, Israel}
\author[L.A. Shuster]{L.A. Shuster}
 \address{Department of Mathematics,
 Bar-Ilan University, 52900 Ramat Gan, Israel}

\begin{abstract}
We consider  an equation
\begin{equation}\label{1}
-(r(x)y'(x))'+q(x)y(x)=f(x),\quad x\in \mathbb R
\end{equation}
where $f\in L_p(\mathbb R),$\ $p\in(1,\infty)$ and
\begin{equation}\label{2}
  r>0, \  q\ge0, \  \frac{1}{r}\in L_1^{\loc}(\mathbb R), \ q\in
 L_1^{\loc}(\mathbb R),
\end{equation}
 \begin{equation}\label{3}
 \int_{-\infty}^0\frac{dt}{r(t)}=\int_0^\infty\frac{dt}{r(t)}=\infty.
 \end{equation}
 By a solution of \eqref{1} we mean any function $y$ which is
 absolutely continuous together with $ry'$ and satisfies \eqref{1}
 almost everywhere on $\mathbb R.$ Under conditions \eqref{2}--\eqref{3},
 we give a criterion for correct solvability of \eqref{1} in
 $L_p(R)$,\ $p\in(1,\infty).$

\end{abstract}

\maketitle

\baselineskip 20pt

\section{Introduction}\label{introduction}
\setcounter{equation}{0} \numberwithin{equation}{section}

In the present paper, we consider an equation
\begin{equation}\label{1.1}
-(r(x)y'(x))'+q(x)y(x)=f(x),\quad x\in \mathbb R
\end{equation}
 where $f\in L_p(\mathbb R),$\ $(L_p(\mathbb R):=L_p),$\ $p\in(1,\infty)$ and
\begin{equation}\label{1.2}
  r>0, \  q\ge0, \  \frac{1}{r}\in L_1^{\loc}(\mathbb R), \ q\in
 L_1^{\loc}(\mathbb R).
\end{equation}

In the sequel, by a solution of equation \eqref{1.1}, we mean any
function $y$ which is absolutely continuous together with $ry'$
and satisfies \eqref{1.1} almost everywhere on $\mathbb R.$ We say
that equation \eqref{1.1} is correctly solvable in a given space
$L_p,$\ $p\in[1,\infty)$ if the following assertions I)--II) hold
(see \cite[Ch.III, \S6, no.2]{1}):
\begin{enumerate}
 \item[1)] for any function $f\in L_p$, there exists a unique
 solution of \eqref{1.1}, $y\in L_p;$

\item[2)] there exists an absolute constant $c(p)\in(0,\infty)$
such that the solution of \eqref{1.1}, $y\in L_p,$ satisfies the
inequality
  \begin{equation}\label{1.3}
   \|y\|_p\le c(p)\|f\|_p,\qquad \forall f\in L_p\ (\|f\|_p:=\|f\|_{L_p}).
 \end{equation}
\end{enumerate}

Our goal is to find exact requirements of $r$ and $q$ which, for a
given $p\in (1,\infty)$, guarantee correct solvability of
\eqref{1.1} in $L_p$. In the sequel, for brevity, this is referred
to as ``problem I)--II)" or ``question on I)--II)". Note that
problem I)--II) can be reformulated in other terms (see
\cite{3,2}).

To this end, let us introduce the set $\mathcal D_p$ and the
operator $\mathcal L_p:$
$$\mathcal D_p=\{y\in L_p: y,ry'\in AC^{\loc}(\mathbb R),\quad
-(ry')'+qy\in L_p\},$$
$$\mathcal L_py=-(ry')'+qy,\quad y\in\mathcal D_p.$$
The linear operator $\mathcal L_p$ is called the maximal
Sturm-Liouville operator, and problem \linebreak I)--II) is,
evidently, equivalent to the problem of existence and boundedness
of the operator $\mathcal L_p^{-1}: L_p\to L_p,$ i.e., to the
problem of continuous invertibility of the operator $\mathcal L_p$
(see \cite{3}). The question on I)--II) in the first or second
formulation was studied in \cite{4,5,6} for $r\equiv 1$ and in
\cite{2,7} for $r\not\equiv1.$ Finally, note that the problem of
continuous invertibility of the minimal Sturm-Liouville operator
$\mathcal L_{o,p}$ was considered in \cite{2,3,8}. This operator
is defined as the closure in $L_p$ of the operator $\mathcal
L_{o,p}':$
$$\mathcal L_{o,p}'=-(ry')'+qy,\quad y\in\mathcal D_{o,p}'$$
where the set $\mathcal D_{o,p}'$ consists of all finitary
functions belonging to $\mathcal D_p.$ The operator $\mathcal
L_{o,p}$ was studied in \cite{2,3} for $p\in[1,\infty),$ and in
\cite{8} for $p=2.$ See \cite{3} for a brief survey of the work on
continuous invertibility of the Sturm-Liouville operators of both
types.

Let us now return to the initial question on I)--II) and present
our results.

\begin{thm}\label{thm1.1}
Suppose that assumptions \eqref{1.2} hold and, in addition,
\begin{equation}\label{1.4}
\int_{-\infty}^0\frac{dt}{r(t)}=\int_0^\infty\frac{dt}{r(t)}=\infty.
\end{equation}
Then if equation \eqref{1.1} is correctly solvable in $L_p,$\
$p\in[1,\infty),$ the following conditions hold:
\begin{equation}\label{1.5}
\int_{-\infty}^x q(t)dt>0,\quad \int_x^\infty q(t)dt>0\quad
\forall x\in \mathbb R,
\end{equation}
\begin{equation}\label{1.6}
\lim_{|d|\to \infty}\int_{x-d}^x\frac{dt}{r(t)}\cdot\int_{x-d}^x
q(t)dt=\infty\quad \forall x\in \mathbb R.
\end{equation}
\end{thm}

\begin{remark}
Conditions \eqref{1.5} and \eqref{1.6} were introduced in \cite{6}
and \cite{7}.
\end{remark}

Below we need the following lemma.

\begin{lem}\label{lem1.3}
\cite{7} \ Under conditions \eqref{1.2} and \eqref{1.5}, the
equation
\begin{equation}\label{1.7}
(r(x)z'(x))'=q(x)z(x),\quad x\in\mathbb R
\end{equation}
has a fundamental system of solutions (FSS) $\{u,v\}$ with the
following properties:
\begin{equation}\label{1.8}
v(x)>0,\ u(x)>0, \ v'(x)\ge0,\ u'(x)\le0\quad \forall x\in\mathbb
R,
\end{equation}
\begin{equation}\label{1.9}
r(x)[v'(x)u(x)-u'(x)v(x)]=1,\quad \forall x\in\mathbb R,
\end{equation}
\begin{equation}\label{1.10}
\lim_{x\to-\infty}\frac{v(x)}{u(x)}=\lim_{x\to+\infty}\frac{u(x)}{v(x)}=0,
\end{equation}
\begin{equation}\label{1.11}
\int_{-\infty}^0\frac{dt}{r(t)v^2(t)}=\int_0^\infty\frac{dt}{r(t)u^2(t)}=\infty,\quad
\int_0^\infty\frac{dt}{r(t)v^2(t)}<\infty,\quad
\int_{-\infty}^0\frac{dt}{r(t)u^2(t)}<\infty.
\end{equation}
Moreover, properties \eqref{1.8}--\eqref{1.11} determine the FSS
$\{u,v\}$ uniquely up to positive constant factors inverse one to
another.
\end{lem}

The FSS from \lemref{lem1.3} will be denoted by $\{u,v\}$ in the
sequel. This FSS will allow us to define the main tools of the
present research --- the Green function $G(x,t)$ and the Green
integral operator $G:$
\begin{equation}\label{1.12}
G(x,t)=\begin{cases} u(x)v(t), &\  x\ge t,\quad
 x,t\in\mathbb R\\
u(t)v(x)& \  x\le t,\quad
 x,t\in\mathbb R\\
\end{cases}
\end{equation}
\begin{equation}\label{1.13}
(Gf)(x)=\int_{-\infty}^\infty G(x,t)f(t)dt,\quad f\in L_p,\quad
x\in\mathbb R.
\end{equation}
\begin{thm}\label{thm1.4}
Suppose that conditions \eqref{1.2} and \eqref {1.5} hold, and let
$p\in[1,\infty).$ Then equation \eqref{1.1} is correctly solvable
in $L_p$ if and only if the operator $G: L_p\to L_p$ is bounded.
\end{thm}
\begin{cor}\label{cor1.5} Suppose that conditions \eqref{1.2} and
\eqref{1.5} hold and equation \eqref{1.1} is correctly solvable in
$L_p$,\ $p\in[1,\infty).$ Then for any function $f\in L_p$ the
solution $y\in L_p$ of \eqref{1.1} is of the form
\begin{equation}\label{1.14}
y(x)=\int_{-\infty}^\infty G(x,t)f(t)dt,\quad x\in\mathbb
R.\end{equation}
\end{cor}

We now need the following assertion.

\begin{lem}\label{lem1.6} \cite{7} Suppose that conditions
\eqref{1.2} and \eqref{1.6} hold. Then for any given $x\in \mathbb
R$ each of the equations
\begin{equation}\label{1.15}
\int_{x-d}^x\frac{dt}{r(t)}\cdot \int_{x-d}^x q(t)dt=1,\quad
\int_x^{x+d}\frac{dt}{r(t)}\cdot\int_x^{x+d}q(t)dt=1
\end{equation}
in $d\ge0$ has a unique finite positive solution. Denote the
solutions of \eqref{1.15} by $d_1(x),$ $d_2(x)$, respectively. For
$x\in\mathbb R$ let us introduce the following auxiliary
functions:
\begin{equation}\label{1.16}
\varphi(x)=\int_{x-d_1(x)}^x\frac{dt}{r(t)},\qquad
\psi(x)=\int_x^{x+d_2(x)}\frac{dt}{r(t)}
\end{equation}
\begin{equation}\label{1.17}
h(x)=\frac{\varphi(x)\psi(x)}{\varphi(x)+\psi(x)}\quad
\left(\equiv\left(\int_{x-d_1(x)}^{x+d_2(x)}q(t)\right)^{-1}\right).
\end{equation}
Furthermore, for every $x\in\mathbb R$, the equation $d\ge0$ in
\begin{equation}\label{1.18}
\int_{x-d}^{x+d}\frac{dt}{r(t)h(t)}=1
\end{equation}
has a unique finite positive function. Denote it by $d(x).$ The
function $d(x)$ is continuous for $x\in\mathbb R.$ In addition,
$$\lim_{x\to-\infty}(x+d(x))=-\infty,\qquad
\lim_{x\to+\infty}(x-d(x))=\infty.$$
\end{lem}

\begin{remark}\label{rem1.7}
Various auxiliary functions similar to the functions in
\lemref{lem1.6} were introduced by M. Otelbaev (see \cite{9}).
\end{remark}

Let us now state the main result of the present paper.

\begin{thm}\label{thm1.8}
Suppose conditions \eqref{1.2} and \eqref{1.4} hold. Then equation
\eqref{1.1} is correctly solvable in $L_p,$\ $p\in(1,\infty)$ if
and only if condition \eqref{1.5} holds and $B<\infty.$ Here
\begin{equation}\label{1.19}
B=\sup_{x\in\mathbb R} h(x)d(x).
\end{equation}

Moreover, this criterion reduces to the unique condition
$B<\infty$ if condition \eqref{1.4} is replaced with condition
\eqref{1.6}. Finally, in any case one of the following assertions
holds:
\begin{enumerate}
\item[$\alpha)$] for every $p\in(1,\infty),$ equation \eqref{1.1}
is correctly solvable in $L_p;$
 \item[$\beta)$] for all
$p\in(1,\infty),$ equation \eqref{1.1} is not correctly solvable
in $L_p.$
\end{enumerate}
\end{thm}

\begin{cor}\label{cor1.9} \cite{7} Suppose conditions \eqref{1.2}
and \eqref{1.6} hold. Then for every $p\in(1,\infty)$, equation
\eqref{1.1} is correctly solvable in $L_p$ if $A>0.$ Here
\begin{equation}\label{1.20}
A=\inf_{x\in\mathbb R}\frac{1}{2d(x)}\int_{x-d(x)}^{x+d(x)}q(t)dt.
\end{equation}
\end{cor}

\begin{cor}\label{cor1.10} \cite{6}
Let $r\equiv 1,$ and suppose that the function $q$ satisfies
condition \eqref{1.2}. Then equation \eqref{1.1} is correctly
solvable in $L_p,$\ $p\in(1,\infty)$ if and only if there exists
$a\in(0,\infty)$ such that $m(a)>0.$ Here
\begin{equation}\label{1.21}
m(a)=\inf_{x\in\mathbb R} \int_{x-a}^{x+a}q(t)dt.
\end{equation}
\end{cor}

\begin{remark}\label{rem1.11}
Under some additional assumptions on the functions $\varphi$ and
$\psi$, \thmref{thm1.8} was obtained in \cite{10}. Corollary
\ref{cor1.9} remains true also for $p=1$ (see \cite{7}). Corollary
\ref{cor1.10} remains true for $p=1$ and $p=\infty$ (here
$L_\infty(\mathbb R):=C(\mathbb R)$, see \cite{6}).
\end{remark}

\begin{remark}\label{rem1.12}
See \S3 for the proofs of all the above statements. In \S4 we give
examples of applications of \thmref{thm1.8} to concrete equations.
\S2 contains a list of various facts used in \S\S 3--4.
\end{remark}

To conclude this introductory section, note that the results and
methods presented here allow us to find exact conditions for:
\begin{enumerate}
\item[1)] correct solvability of equation \eqref{1.1} in the
spaces $L_1(\mathbb R)$ and $C(\mathbb R)$; \item[2)] compactness
of the operator $\mathcal L_p^{-1}: L_p\to L_p$ for
$p\in[1,\infty);$ \item[3)] separability of equation \eqref{1.1}
in $L_p$,\ $p\in(1,\infty)$ (the problem of Everitt-Giertz, see
\cite{14,12,13}; and see \cite{2,7} for the case $p=1).$ The
solutions of these problems will appear in our forthcoming papers.
\end{enumerate}

\section{Preliminaries}

\begin{thm}\label{thm2.1} \cite{16} Under conditions \eqref{1.2}
and \eqref{1.5} on the FSS $\{u,v\}$ of equation \eqref{1.7} and
the Green function $G(x,t)$ (see \eqref{1.12}), the Davies-Harrell
representations hold:
\begin{equation}
\begin{aligned}\label{2.1}
&u(x)=\sqrt{\rho(x)}\exp\left(-\frac{1}{2}\int_{x_0}^x\frac{dt}{r(t)\rho(t)}\right),\quad
x\in\mathbb R
\\
&v(x)=\sqrt{\rho(x)}\exp\left(\frac{1}{2}\int_{x_0}^x\frac{dt}{r(t)\rho(t)}\right),\quad
x\in\mathbb R
\end{aligned}
\end{equation}

\begin{equation}\label{2.3}
G(x,t)=\sqrt{\rho(x)\rho(t)}\exp\left(-\frac{1}{2}\left|\int_{t}^x\frac{d\xi}{r(\xi)\rho(\xi)
}\right|\right).
\end{equation}
\end{thm}

Here $x,t\in\mathbb R,$\ $x_0$ is the unique solution of the
equation $u(x)=v(x)$ in $\mathbb R$ and
\begin{equation}\label{2.4}
\rho(x)=u(x)v(x),\quad x\in\mathbb R.
\end{equation}

\begin{remark}\label{rem2.2}
Representations \eqref{2.1} and \eqref{2.3} were obtained in
\cite{16} for $r\equiv1.$ \thmref{thm2.1} was proved in \cite{7}.
\end{remark}

\begin{lem}\label{lem2.3} \cite{7} Suppose conditions \eqref{1.2}
and \eqref{1.6} hold. Then for $x\in\mathbb R$ we have the
following inequalities (see \eqref{1.7} and \eqref{1.16}:
\begin{equation}
\begin{aligned}\label{2.5}
2^{-1}v(x)&\le r(x)v'(x)\varphi(x)\le 2v(x), \\
2^{-1}u(x)&\le r(x)|u'(x)|\psi(x)\le 2u(x).
\end{aligned}
\end{equation}
\end{lem}

\begin{remark}\label{rem2.4}
Two-sided, sharp by order, a priori estimates of type \eqref{2.5}
first appeared in \cite{4} (for $r\equiv1$ and some additional
requirements to $q$). Under the conditions \eqref{1.2} and
$\inf_{x\in\mathbb R} q(x)>0,$ estimates similar to \eqref{2.5},
with other, more complicated, auxiliary functions, were given in
\cite{2}.
\end{remark}

\begin{lem}\label{lem2.5} \cite{7} Suppose conditions \eqref{1.1}
and \eqref{1.6} hold. The for $x\in\mathbb R$, we have the
inequalities (see \eqref{1.17} and \eqref{2.4}:
\begin{equation}\label{2.6}
2^{-1}h(x)\le\rho(x)\le 2h(x).
\end{equation}
\end{lem}

\begin{lem}\label{lem2.6} \cite{7} Let $r\equiv1$ and suppose that
$q$ satisfies conditions \eqref{1.1} and \eqref{1.5}. Then for
every given $x\in\mathbb R$ the equation
\begin{equation}\label{2.7}
d\int_{x-d}^{x+d} q(t)dt=2
\end{equation}
in $d\ge0$ has a unique finite positive solution. Denote it by
$\tilde d(x).$ We have the inequalities (see \eqref{2.4}:
\begin{equation}\label{2.8}
4^{-1}\tilde d(x)\le \rho(x)\le 3\cdot 2^{-1}\tilde d(x),\quad
x\in\mathbb R.
\end{equation}
\end{lem}

\begin{remark}\label{rem2.7}
Two-sided, sharp by order, a priori estimates of the function
$\rho$ first appeared in \cite{15} (under some additional
requirements to $r$ and $q).$ Therefore, we call all inequalities
of such a form Otelbaev inequalities. Note that in \cite{15} there
were used other, more complicated auxiliary functions than $h(x)$
and $\tilde d(x).$ The function $\tilde d(x)$ was introduced by M.
Otelbaev (see \cite{9}).
\end{remark}

\begin{lem}\label{lem2.8} \cite{7}
Suppose conditions \eqref{1.2} and \eqref{1.6} hold. Then for all
$x \in\mathbb R$ and $t\in[x-d(x),x+d(x)]$ (see \eqref{1.18}), we
have the inequalities (see \eqref{2.4}):
 \begin{equation}\label{2.9} e^{-2}  \rho(x)\le \rho(t)\le
 e^2\rho(x).
 \end{equation}
\end{lem}

\begin{defn}\label{defn2.9} \cite{7} We say that a system of
segments $\{\Delta_n\}_{n\in\mathbb N'},$\ $\mathbb
N'=\{\pm1,\pm2,\dots\}$ forms an $\mathbb R(x)$-covering of
$\mathbb R$ if the following assertions hold (see \eqref{1.18}):
\begin{enumerate}
\item[1)] $\Delta_n=[\Delta_n^{-},\Delta_n^+]\doe
[x_n-d(x_n),x_n+d(x_n)],\quad n\in\mathbb N';$
 \item[2)]
$\Delta_{n+1}^-=\Delta_n^+,\quad \text{if}\quad n\ge1;\
\Delta_{n-1}^+=\Delta_n^-\quad\text{if}\quad n\le-1;$
 \item[3)]
$\Delta_{1}^-=\Delta_{-1}^+=x,\quad
\bigcup\limits_{n\ne0}\Delta_n=\mathbb R.$
\end{enumerate}
\end{defn}

\begin{lem}\label{lem2.10} \cite{7} Under conditions \eqref{1.2}
and \eqref{1.6}, for every $x\in\mathbb R$, there exists an
$\mathbb R(x)$-covering of $\mathbb R.$
\end{lem}

\begin{remark}\label{rem2.11} Assertions similar to
\lemref{lem2.10} were introduced and systematically studied by M.
Otelbaev (see \cite{9}).
\end{remark}

\begin{lem}\label{lem2.12} \cite{7} Suppose conditions \eqref{1.2}
and \eqref{1.5} hold. Then equation \eqref{1.7} has no solutions
$z\in L_p,$ \ $p\in [1,\infty)$ apart from $z\equiv 0.$
\end{lem}

\begin{thm}\label{thm2.13} Let $\mu,\theta$ be continuous positive
functions in $\mathbb R,$ and let $\mathcal K$ be an integral
operator:
\begin{equation}\label{2.10}
(\mathcal K f)(t)=\mu(t)\int_t^\infty\theta(\xi)f(\xi)d\xi,\quad
t\in\mathbb R.
\end{equation}
Then for $p\in(1,\infty)$ the operator $\mathcal K: L_p\to L_p$ is
bounded if and only if $H_p<\infty.$ Here
$H_p=\sup\limits_{x\in\mathbb R}H_p(x),$
\begin{equation}\label{2.11}
H_p(x)=\left[\int_{-\infty}^x \mu(t)^pdt\right]^{1/p}\left[
\int_x^\infty\theta(t)^{p'}dt\right]^{1/p'},\quad
p'=\frac{p}{p-1},\quad x\in\mathbb R.
\end{equation}

Moreover, the following inequalities hold:
\begin{equation}\label{2.12}
H_p\le\|\mathcal K\|_{p\to p}\le (p)^{1/p}(p')^{1/p'}H_p.
\end{equation}
\end{thm}

\begin{thm}\label{thm2.14}
Let $\mu,\theta$ be continuous positive functions in $\mathbb R,$
and let $\tilde{\mathcal K}$ be an integral operator:
\begin{equation}\label{2.13}
 (\tilde{\mathcal
K}f)(t)=\mu(t)\int_{-\infty}^t\theta(\xi)f(\xi)d\xi,\quad
t\in\mathbb R.
\end{equation}
Then for $p\in(1,\infty)$ the operator $\tilde{\mathcal K}: L_p\to
L_p$ is bounded if and only if $\tilde H_p<\infty.$ Here $\tilde
H_p=\sup\limits_{x\in\mathbb R}\tilde H_p(x),$
\begin{equation}\label{2.14}
\tilde
H_p(x)=\left[\int_{-\infty}^x\theta(t)^{p'}dt\right]^{1/p'}\left[\int_x^\infty\mu(t)^pdt\right]^{1/p},\quad
p'=\frac{p}{p-1},\quad x\in\mathbb R.
\end{equation}

Moreover, the following inequalities hold:
\begin{equation}\label{2.15}
\tilde H_p\le\|\tilde{\mathcal K}\|_{p\to p}\le
(p)^{1/p}(p')^{1/p'}\tilde H_p.
\end{equation}
\end{thm}

\begin{remark}\label{rem2.15} Theorems \ref{thm2.13} and
\ref{thm2.14} follow from a  Hardy type inequality (see
\cite{11}). In particular, see \cite{17} for such a proof. See
\cite{18} for the original direct proof of these theorems (under
weaker requirements of $\mu$ and $\theta).$
\end{remark}

\begin{thm}\label{thm2.16} \cite[Ch.V, \S2, no.5]{19} Let
$-\infty\le a<b\le \infty,$ let $\mathcal K(s,t)$ be a continuous
function for $s,t\in[a,b]$, and let $\mathcal K$ be an integral
operator
\begin{equation}\label{2.16}
(\mathcal K f)(t)=\int_a^b\mathcal K(s,t)f(s)ds,\qquad t\in [a,b].
\end{equation}
Then
\begin{equation}\label{2.17} \|\mathcal K\|_{L_1(a,b)\to L_1(a,b)}=\sup_{s\in[a,b]} \int_a^b|\mathcal K
(s,t)|dt.
\end{equation}
\end{thm}

\section{Proofs}

\begin{proof}[Proof of \thmref{thm1.1}]

Assume the contrary. Then for some $p\in[1,\infty)$, equation
\eqref{1.1} is correctly solvable in $L_p,$ and there exists a
point $x_0\in\mathbb R$ such that
\begin{equation}\label{3.1}
\int_{x_0}^\infty q(t)dt=0\quad\Rightarrow\quad
q(t)=0\quad\text{almost everywhere on} \quad (x_0,\infty).
\end{equation}
Let $x_1\gg\max\{1,|x_0|\},$ and let
\begin{equation}\label{3.2}
f(x)=\begin{cases} -1,&\quad\text{if}\quad x\in [x_0,x_1)\\
0,&\quad\text{if}\quad x\notin [x_0,x_1)\end{cases}
\quad\Rightarrow\quad \|f\|_p=(x_1-x_0)^{1/p}.
\end{equation}
By \eqref{3.1} and \eqref{3.2}, for $x\ge x_1,$ equation
\eqref{1.1} is of the form
\begin{align}
-(r(x)y'(x))'=0,\quad &x\ge x_1\qquad \Rightarrow\label{3.3}\\
y(x)=c_1+c_2\int_{x_1}^x\frac{dt}{r(t)},\quad& x\ge
x_1.\label{3.4}
\end{align}
Here $c_1,c_2$ are some constants. From I)--II) (see \S1) it
follows that
\begin{equation}\label{3.5}
\|y\|_{L_p(x_1,\infty)}\le\|y\|_p\le c(p)\|f\|_p<\infty.
\end{equation}
Then, according to \eqref{3.5} and \eqref{1.4}, we conclude that
$c_1=c_2=0.$ Indeed, if $c_2\ne0,$ then for all $x\gg x_1$ we get
$$|y(x)|\ge
|c_2|\int_{x_1}^x\frac{dt}{r(t)}\left[1-\left|\frac{c_1}{c_2}\right|\left(\int_{x_1}^x\frac{dt}{r(t)}
\right)^{-1}\right]\ge\frac{|c_2|}{2}>0$$ and therefore
$\|y\|_p=\infty.$ Contradiction. Hence $c_2=0.$ But then also
$c_1=0$ because otherwise $y=c_1\notin L_p.$ Thus $y(x)\equiv0$
for $x\ge x_1.$ Hence in the case \eqref{3.2}, equation
\eqref{1.1} has solution $y\in L_p$ which, for $x\in [x_0,x_1],$
satisfies the relations
\begin{equation}\label{3.6}
-(r(x)y'(x))'=-1,\ y(x_1)=y'(x_1)=0\quad\Rightarrow \quad
y(x)=\int_x^{x_1}\frac{x_1-t}{r(t)}dt,\quad x\in [x_0,x_1].
\end{equation}

Let us first consider the case $p=1.$ From \eqref{3.6} it follows
that
\begin{align}
\|y\|_{L_1(x_0,x_1)}&=\int_{x_0}^{x_1}\left(\int_x^{x_1}\frac{x_1-t}{r(t)}dt\right)dx
=(x-x_0)\int_x^{x_1}\frac{x_1-t}{r(t)}dt\bigm|_{x_0}^{x_1}+\int_{x_0}^{x_1}\frac{(x-x_0)(x_1-x)}{
 r(x)}dx\nonumber\\
&= \int_{x_0}^{x_1}\frac{(x-x_0)(x_1-x)}{r(x)}dx.\label{3.7}
\end{align}
Since $x_1\gg\max\{1,|x_0 |\},$ according to \eqref{3.7} we get
\begin{align}
\|y\|_{L_1(x_0,x_1)}&\ge
\int_{x_0+1}^{x_1-1}\frac{(x-x_0)(x_1-x)}{r(x)}dx\ge(x_1-x_0-1)
\int_{x_0+1}^{x_1-1}
\frac{dx}{r(x)}\nonumber\\
&\ge\frac{x_1-x_0}{2}\int_{x_0+1}^{x_1-1}\frac{dx}{r(x)}.\label{3.8}
\end{align}
But then from \eqref{3.8} and \eqref{1.3}, it follows that
\begin{gather}
\frac{x_1-x_0}{2}\int_{x_0+1}^{x_1-1}\frac{dx}{r(x)}\le\|y\|_{L_1(x_0,x_1)}\le\|y\|_1\le
c(1)\|f\|_1=c(1)(x_1-x_0)\quad\Rightarrow\nonumber\\
\int_{x_0+1}^{x_1-1}\frac{dx}{r(x)}\le 2c(1)<\infty.\label{3.9}
\end{gather}

Since $x_1$ can be chosen arbitrarily large, \eqref{3.9}
contradicts \eqref{1.4}, and in the case $p=1$ the theorem is
proven.

Consider the case $p\in (1,\infty).$ From \eqref{3.6}, for
$x\in[x_0,x_0+1],$ it follows that
\begin{align}
y(x)&=\int_x^{x_0+1}\frac{x_1-t}{r(t)}dt+\int_{x_0+1}^{x_1}\frac{x_1-t}{r(t)}dt
>\int_{x_0+1}^{x_1}\frac
{x_1-t}{r(t)}dt>\int_{x_0+1}^{x_0+2}\frac{x_1-t}{r(t)}dt\nonumber\\
&\ge\int_{x_0+1}^{x_0+2}\frac{x_1-x_0-2}{r(t)}dt\ge\frac{x_1-x_0}{2}
\int_{x_0+1}^{x_0+2}\frac{dt}{r(t)}\quad
\Rightarrow\nonumber\\
&\qquad \|y\|_{L_p(x_0,x_0+1)}\ge
\frac{x_1-x_0}{2}\int_{x_0+1}^{x_0+2}\frac{dt}{r(t)}.\label{3.10}
\end{align}
Then, according to \eqref{3.10} and \eqref{1.3}, we get
\begin{align}
\frac{x_1-x_0}{2}\int_{x_0+1}^{x_0+2}\frac{dt}{r(t)}&\le\|y\|_{L_p(x_0,x_0+1)}\le\|y\|_p\le
c(p)\|f\|_p=c(p)(x_1-x_0)^{1/p}\nonumber\\
&\Rightarrow
(x_1-x_0)^{1/p'}\int_{x_0+1}^{x_0+2}\frac{dt}{r(t)}\le
2c(p)<\infty,\quad p'=\frac{p}{p-1}.\label{3.11}
\end{align}
Since $x_1$ can be taken arbitrarily large, \eqref{3.11}
contradicts \eqref{1.2}. Thus inequalities \eqref{1.5} are proven.
It remains to notice that \eqref{1.6} follows from \eqref{1.5} and
\eqref{1.4}.
\end{proof}

\renewcommand{\qedsymbol}{}
\begin{proof}[Proof of \thmref{thm1.4}]  Necessity.

We need the following (maybe commonly known) assertion.
\end{proof}

\begin{lem}\label{lem3.1}
Suppose conditions \eqref{1.2} and \eqref{1.5} hold. Consider the
integral equations
\begin{equation}\label{3.12}
(G_1f)(x)=u(x)\int_{-\infty}^x v(t)f(t)dt,\quad x\in\mathbb R;
\end{equation}
\begin{equation}\label{3.13}
(G_2f)(x)=v(x)\int_{x}^\infty u(t)f(t)dt,\quad x\in\mathbb R;
\end{equation}
For $p\in[1,\infty)$ the following inequalities hold:
\begin{equation}\label{3.14}
\frac{\|G_1\|_{p\to p}+\|G_2\|_{p\to p}}{2}\le\|G\|_{p\to p}\le
\|G_1\|_{p\to p}+\|G_2\|_{p\to p}.
\end{equation}
\end{lem}

\renewcommand{\qedsymbol}{\openbox}

\begin{proof} The upper estimate for $\|G\|_{p\to p}$ follows
from the triangle inequality for norms. Furthermore, the following
relations are obvious:
\begin{align*}
\|G_1f\|_p&=\left[\int_{-\infty}^\infty
u(x)^p\left|\int_{-\infty}^x v(t)f(t)dt\right|^pdx\right]\\
&\le \left[\int_{-\infty}^\infty u(x)^p\left(\int_{-\infty}^x
v(t)|f(t)|dt\right)^pdx\right]^{1/p}\\
&\le\left\{\int_{-\infty}^\infty\left[u(x)\int_{-\infty}^x
v(t)|f(t)|dt+v(x)\int_x^\infty
u(t)|f(t)|dt\right]^pdx\right\}^{1/p}\\
&=\|G|f|\,\|_p\le \|G\|_{p\to p}\cdot \|f\|_p\quad\Rightarrow\quad
\|G_1\|_{p\to p}\le \|G\|_{p\to p}.
\end{align*}

In a similar way, one can check that $\|G_2\|_{p\to p}\le
\|G\|_{p\to p}.$ These inequalities imply the lower estimate of
\eqref{3.14}.
\end{proof}

Let us now go to the proof of the theorem. Suppose that equation
\eqref{1.1} is correctly solvable in $L_p$ for some
$p\in[1,\infty).$ Take an arbitrary pair of numbers $x_1,x_2$
$(x_1\le x_2),$ and for any function $f\in L_p$ set
$$\tilde f(x)=\begin{cases} f(x),&\quad\text{if}\quad
x\in[x_1,x_2]\\
0,&\quad\text{if}\quad x\notin[x_1,x_2]\end{cases}$$

Let $\tilde y$ be a solution of \eqref{1.1} with right-hand side
$\tilde f$ such that $\tilde y\in L_p.$ It is easy to see that
$\tilde y$ is of the form
\begin{equation}
\begin{aligned}\label{3.15}
\tilde y(x)&=c_1u(x)+c_2v(x)+u(x)\int_{-\infty}^xv(t)\tilde
f(t)dt\\
&\quad + v(x)\int_x^\infty u(t)\tilde f(t)dt,\quad x\in\mathbb R.
\end{aligned}
\end{equation}
Here $c_1,c_2$ are some constants. Indeed, if the integrals in
\eqref{3.15} exist, then representation \eqref{3.15} holds true
for $\tilde y $ because this formula represents the general
solution of \eqref{1.1}. The existence of these integrals follows
from the definition of $\tilde f.$ For example, for
$p\in(1,\infty) $ and $x\in\mathbb R,$ using H\"older's inequality
and the definition of $\tilde f,$ we obtain
\begin{equation}
\begin{aligned}\label{3.16}
&\int_{-\iy}^xv(t)|\tilde f(t)|dt \le\int_{x_1}^{x_2}v(t)|\tilde
f(t)|dt
 \le\left(\int_{x_1}^{x_2}v(t)^{p'}dt\right)^{1/p'} \|f\|_p,\quad
p'=\frac{p}{p-1}; \\
&\int_{-\iy}^xu(t)|\tilde f(t)|dt \le\int_{x_1}^{x_2}u(t)|\tilde
f(t)|dt
 \le\left(\int_{x_1}^{x_2}u(t)^{p'}dt\right)^{1/p'} \|f\|_p,\quad
p'=\frac{p}{p-1}.
\end{aligned}
\end{equation}
Thus, the integrals exist and formula \eqref{3.15} holds true.

Let us now prove that $c_1=c_2=0.$ Assume the contrary. Let, say,
$c_2\ne0.$ Then for $x\ge x_2$, we have
\begin{align*}
|\tilde y(x)|&\ge
|c_2|v(x)-|c_1|u(x)-u(x)\int_{x_1}^{x_2}v(t)|\tilde f(t)|dt\\
&=|c_2|v(x)\left[1-\left|\frac{c_1}{c_2}\right|\frac{u(x)}{v(x)}-\frac{1}{|c_2|}\
\frac{u(x)}{v(x)}\int_{x_1}^{x_2}v(t)|\tilde f(t)|dt\right].
\end{align*}
The latter inequality, together with \eqref{3.16} and
\eqref{1.10}, implies
$$|\tilde y(x)|\ge \frac{|c_2|}{2}v(x),\quad
x\gg\max\{1,|x_2|\}\Rightarrow \|\tilde y\|_p=\infty,$$ a
contradiction. Hence $c_2=0$ and, similarly, $c_1=0,$ and
therefore
\begin{equation}\label{3.17}
\tilde y(x)=\int_{-\infty}^\infty G(x,t)\tilde
f(t)dt=\int_{x_1}^{x_2}G(x,t)\tilde f(t)dt,\quad x\in \mathbb R.
\end{equation}

Below we consider the cases $p=1$ and $p\in(1,\infty) $
separately. For $p=1,$ using \eqref{3.17} and \eqref{1.3}, we get
$$\|\tilde y\|_{L_1(x_1,x_2)}\le\|\tilde y\|_1\le c(1)\|\tilde
f\|_1=c(1)\|\tilde f\|_{L_1(x_1,x_2)}.$$ Hence the operator $G:
L_1(x_1,x_2)\to L_1(x_1,x_2) $ (see \eqref{3.17}) is bounded, and
its norm does not exceed $c(1)$ (see \thmref{thm2.16}), i.e.:
\begin{equation}\label{3.18}
\|G\|_{L_1(x_1,x_2)\to
L_1(x,x_2)}=\sup_{t\in[x_1,x_2]}\int_{x_1}^{x_2}G(x,t)dx\le c(1).
\end{equation}
Since $x_1,x_2$ are arbitrary numbers, \eqref{3.18} implies the
inequality
$$\sup_{t\in\mathbb R}\int_{-\infty}^\infty G(x,t)dx\le c(1).$$
Then, using \thmref{thm2.16} once again, we get
$$\|G\|_{1\to 1}=\sup_{t\in\mathbb R}\int_{-\infty}^\infty
G(x,t)dx\le c(1)<\infty.$$ Q.E.D.

In the case $p\in(1,\infty),$ set
\begin{gather}
\tilde f(x)=\begin{cases} u(x)^{p'-1}&\quad\text{if}\quad x\in
[x_1,x_2],\quad p'=\frac{p}{p-1}\\ 0&\quad\text{if}\quad
x\notin[x_1,x_2] \end{cases}\quad\Rightarrow\nonumber\\
\|\tilde
f\|_p=\left(\int_{x_1}^{x_2}u(x)^{p(p'-1)}dx\right)^{1/p}=\left(\int_{x_1}^{x_2}u(x)^{p'}dx
\right)^{1/p}<\infty.\label{3.19}
\end{gather}

Furthermore, from \eqref{3.17} it follows that
\begin{align*}
\tilde y(x)&=\int_{-\infty}^\infty G(x,t)\tilde
f(t)dt=u(x)\int_{-\infty}^x v(t)\tilde
f(t)dt+v(x)\int_x^\infty u(t)\tilde f(t)dt\\
&\ge v(x)\int_x^\infty u(t)\tilde f(t)dt\quad\Rightarrow\\
\|\tilde y\|_p^p&\ge \int_{-\infty}^\infty
v(x)^p\left[\int_x^\infty u(t)\tilde f(t)dt\right]^pdx\ge
\int_{-\infty}^{x_1}v(x)^p\left[\int_x^\infty u(t)\tilde
f(t)dt\right]^pdx\\
&\ge\left(\int_{-\infty}^{x_1}v(x)^pdx\right)\left(\int_{x_1}^\infty
u(t)\tilde
f(t)dt\right)^p=\left(\int_{-\infty}^{x_1}v(x)^pdx\right)\left(\int_{x_1}^{x_2}u(t)^{p'}dt\right)^p\quad
\Rightarrow \end{align*}
 \begin{equation}\label{3.20}
\|\tilde y\|_p\ge\left
[\int_{-\infty}^{x_1}v(t)^pdt\right]^{1/p}\left[\int_{x_1}^{x_2}u(t)^{p'}dt\right].
\end{equation}

{}From \eqref{1.3}, \eqref{3.19} and \eqref{3.20}, we now get
\begin{align*}
\left[\int_{-\infty}^{x_1}v(t)^pdt\right]^{1/p}\left[\int_{x_1}^{x_2}u(t)^{p'}dt\right]&\le
\|\tilde y\|_p\le c(p)\|\tilde f\|_p\\
&=c(p)\left[\int_{x_1}^{x_2}u(t)^{p'}dt\right]^{1/p}\quad\Rightarrow
\end{align*}
\begin{equation}\label{3.21}
\left[\int_{-\infty}^{x_1}v(t)^pdt\right]^{1/p}\left[\int_{x_1}^{x_2}u(t)^{p'}dt\right]^{1/p'}\le
c(p)<\infty.
\end{equation}

Since in \eqref{3.21} the numbers $x_1,x_2$ are arbitrary,
\eqref{3.21} implies
$$\sup_{x\in\mathbb R}\left[\int_{-\infty}^x
v(t)^pdt\right]^{1/p}\left[\int_x^\infty
u(t)^{p'}dt\right]^{1/p'}\le c(p)<\infty.$$ By \thmref{thm2.13},
this implies that the operator $G_2: L_p\to L_p$ (see
\eqref{3.13}) is bounded. Similarly, using \thmref{thm2.14}, we
conclude that the operator $G_1: L_p\to L_p$ (see \eqref{3.12}) is
bounded, too. Hence the operator $G: L_p\to L_p$ is bounded in
view of \lemref{lem3.1}.

\begin{proof}[Proof of \thmref{thm1.4}]  Sufficiency.

We shall only prove the theorem in the case $p\in(1,\infty)$ since
for $p=1$ the argument is simpler and goes along similar lines. So
let $\|G\|_{p\to p}<\infty.$ Then $\|G_1\|_{p\to p},$\
$\|G_2\|_{p\to p}<\infty$ by \lemref{lem3.1}, and therefore from
Theorems \ref{thm2.13} and \ref{thm2.14} we obtain the
inequalities:
$$\int_{-\infty}^x v(t)^{p'}dt<\infty,\qquad \int_x^\infty
u(x)^{p'}<\infty\qquad \forall x\in\mathbb R.$$ Together with
H\"older's inequality this implies that for any $x\in\mathbb R$
and any function $f\in L_p$, there exist the integrals
$$\int_{-\infty}^x v(t)f(t)dt,\qquad \int_x^\infty u(t)f(t)dt.$$
Thus for $x\in\mathbb R $ the following function is well-defined:
$$y(x)=u(x)\int_{-\infty}^x v(t)f(t)dt+v(x)\int_x^\infty
u(t)f(t)dt\equiv (Gf)(x).$$ A straightforward computation shows
that $y=Gf$ is a solution of \eqref{1.1} for which inequality
\eqref{1.3} holds because $\|G\|_{p\to p}<\infty.$ To complete the
proof of the theorem, it remains to apply \lemref{lem2.12}.
\end{proof}

\begin{proof}[Proof of Corollary \ref{cor1.5}]
If equation \eqref{1.1} is correctly solvable in $L_p,$\ $p\in
[1,\infty),$ then the operator $G:L_p\to L_p$ is bounded in view
of \thmref{thm1.4}. At this point, to prove the corollary, it is
enough to repeat the argument from the ``sufficiency part" of
\thmref{thm1.4}.
\end{proof}

\renewcommand{\qedsymbol}{}
\begin{proof}[Proof of \thmref{thm1.8}]  Necessity.

Below we need the following simple fact.
\end{proof}

\begin{lem}\label{lem3.2} Suppose we are under the conditions of
\thmref{thm1.8}. Consider the functions
\begin{equation}\label{3.22}
\Phi_1(x)=\left[\int_{-\infty}^x
v(t)^{p'}dt\right]^{1/p'}\cdot\left[\int_x^\infty
u(t)^pdt\right]^{1/p},\quad x\in\mathbb R;
\end{equation}
\begin{equation}\label{3.23}
\Phi_2(x)=\left[\int_{-\infty}^x v(t)^{p
}dt\right]^{1/p}\cdot\left[\int_x^\infty
u(t)^{p'}dt\right]^{1/p'},\quad x\in\mathbb R,
\end{equation}
where $p\in(1,\infty),$\ $p'=\frac{p}{p-1}.$ Then for $x\in\mathbb
R,$ we have
\begin{align}
\Phi_1(x)&=\left[\int_{-\infty}^x
\rho(t)^{p'/2}\exp\left(-\frac{p'}{2}\int_t^x\frac{d\xi}{r(\xi)\rho(\xi)}\right)dt\right]^{1/p'}\nonumber\\
  &\quad \cdot
 \left[
\int_x^{\infty}\rho(t)^{p/2}\exp
 \left(-\frac{p}{2}\int_x^t\frac{d\xi}{r(\xi)\rho(\xi)}\right)dt
 \right]^{1/p}\label{3.24}
\end{align}
\begin{align}
\Phi_2(x)&=\left[\int_{-\infty}^x
\rho(t)^{p/2}\exp\left(-\frac{p}{2}\int_t^x\frac{d\xi}{r(\xi)\rho(\xi)}\right)dt\right]^
{1/p} \nonumber\\
  &\quad \cdot
 \left[
\int_x^{\infty}\rho(t)^{p'/2}\exp
 \left(-\frac{p'}{2}\int_x^t\frac{d\xi}{r(\xi)\rho(\xi)}\right)dt
 \right]^{1/p'}.\label{3.25}
\end{align}
\end{lem}

\renewcommand{\qedsymbol}{\openbox}
\begin{proof}
Equalities \eqref{3.24} and \eqref{3.25} are proved in the same
way. For example, let us prove \eqref{3.25}. Let us substitute
\eqref{2.1} into \eqref{3.23}:
\begin{align*}
\Phi_2(x)&=\left[\int_{-\infty}^x\rho(t)^{p/2}\exp\left(\frac{p}{2}\int_{x_0}^t\frac{d\xi}{r(\xi)\rho(\xi)}\right)
dt\right]^{1/p}\left[\int_x^\infty\rho(t)^{p'/2}\exp\left(-\frac{p'}{2}\int_{x_0}^t\frac{d\xi}{r(\xi)\rho(\xi)}
\right)dt\right]^{1/p'}\\
&=\left[\int_{-\infty}^x\rho(t)^{p/2}\exp\left(-\frac{p}{2}\int_t^x\frac{d\xi}{r(\xi)\rho(\xi)}\right)\cdot
\exp\left(\frac{p}{2}\int_{x_0}^x\frac{d\xi}{r(\xi)\rho(\xi)}\right)dt\right]^{1/p}\\
&\quad \cdot \left[\int_x^\infty
\rho(t)^{p'/2}\exp\left(-\frac{p'}{2}\int_x^t\frac{d\xi}{r(\xi)\rho(\xi)}\right)\exp\left(-\frac{p'}{2}
\int_{x_0}^x\frac{d\xi}{r(\xi)\rho(\xi)}\right)dt\right]^{1/p'}\\
&=\left[\int_{-\infty}^x\rho(t)^{p/2}\exp\left(-\frac{p}{2}\int_t^x\frac{d\xi}{r(\xi)\rho(\xi)}\right)dt\right]
^{1/p}\\
&\quad
\cdot\left[\int_x^\infty\rho(t)^{p'/2}\exp\left(-\frac{p'}{2}\int_x^t\frac{d\xi}{r(\xi)\rho(\xi)}\right)
dt\right]^{1/p'}.
\end{align*}
\end{proof}

Let us now go to the proof of the theorem. Suppose conditions
\eqref{1.2} and \eqref{1.4} hold, and for some $p\in(1,\infty)$
equation \eqref{1.1} is correctly solvable in $L_p.$ Then by
\thmref{thm1.1}, relations \eqref{1.5} and \eqref{1.6} hold. Hence
all auxiliary functions from \lemref{lem1.6} are well-defined.
Furthermore, by \thmref{thm1.4} the operator $G: L_p\to L_p$ is
bounded, and therefore the operator $G_2: L_p\to L_p$ (see
\lemref{lem3.1}) is bounded, too. In the following relations, we
consecutively use \eqref{3.14}, \eqref{2.12}, \eqref{3.25},
\lemref{lem1.6}, \eqref{2.9}, \eqref{2.6} and \eqref{1.18}:


\begin{align*}
\infty&>2\|G\|_{p\to p}\ge\|G_2\|_{p\to p}\ge \sup_{x\in\mathbb
R}\left[\int_{-\infty}^x v(t)^pdt\right]^{1/p}\left[\int_x^\infty
u(t)^{p'}dt\right]^{1/p'}\\
 &=\sup_{x\in\mathbb R}
 \left[\int_{-\infty}^x\rho(t)^{p/2}\exp\left(-\frac{p}{2}\int_t^x\frac{d\xi}{r(\xi)\rho(\xi)}
 \right)dt\right]^{1/p}\\
 &\quad\cdot
 \left[\int_x^\infty\rho(t)^{p'/2}\exp\left(-\frac{p'}{2}\int_x^t\frac{d\xi}{r(\xi)\rho(\xi)}\right)dt\right]
 ^{1/p'} \\ 
 &\ge\sup_{x\in\mathbb
 R}\left[\int_{x-d(x)}^x\rho(t)^{p/2}\exp\left(-\frac{p}{2}\int_t^x\frac{d\xi}{r(\xi)\rho(\xi)}\right)dt\right]
 ^{1/p}\\
 &\quad\cdot\left[\int_x^{x+d(x)}\rho(t)^{p'/2}\exp\left(-\frac{p'}{2}\int_{x}^t\frac{d\xi}{r(\xi)\rho(\xi)}\right)
 dt\right]^{1/p'}\\
 &\ge \sup_{x\in\mathbb R}
 e^{-4}\rho(x)\exp\left(-\frac{1}{2}\int_{x-d(x)}^{x+d(x)}\frac{d\xi}{r(\xi)\rho(\xi)}
 \right)d(x)^{\frac{1}{p}
 +\frac{1}{p'}}
\\
&\ge 2^{-1}e^{-5}\sup_{x\in\mathbb R}h(x)d(x)=2^{-1}e^{-5}B.
\end{align*}

\begin{proof}[Proof of \thmref{thm1.8}]  Sufficiency.

Suppose conditions \eqref{1.2}, \eqref{1.4}, \eqref{1.5} hold, and
let $p\in(1,\infty).$ Then equality \eqref{1.6} evidently holds,
and therefore all the auxiliary functions from \lemref{lem1.6} are
well-defined. Let us show that the condition $B<\infty$ guarantees
that the operator $G: L_p\to L_p$ is bounded and thus complete the
proof (see \thmref{thm1.4}). By \lemref{lem3.1}, we have
$\|G\|_{p\to p}<\infty$ if $\|G_1\|_{p\to p}<\infty$ and
$\|G_2\|_{p\to p}<\infty.$ The latter inequalities follow from the
estimates
\begin{equation}\label{3.26}
\|G_1\|_{p\to p}\le c(p)B,\qquad \|G_2\|_{p\to p}\le cB.
\end{equation}
Both inequalities in \eqref{3.26} are checked in the same way.
Consider, say, the second one. Let first $p\in(1,2],$ \ $f\in
L_p,$ and denote by $\alpha$ some number from the interval $[0,1)$
which will be chosen later. Then by \lemref{lem1.3}, for any
$x\in\mathbb R,$ we have
\begin{align}
|(G_2f)(x)|&=v(x)\left|\int_x^\infty u(t)f(t)dt\right|\le
v(x)\int_x^\infty u(t)|f(t)|dt\nonumber\\
&\le (v(x)u(x))^\alpha v(x)^{1-\alpha}  \int_x^\infty
u(t)^{1-\alpha}|f(t)|dt\nonumber\\
&=\rho(x)^\alpha v(x)^{1-\alpha}\int_x^\infty
u(t)^{1-\alpha}|f(t)|dt,\quad x\in\mathbb R.\label{3.27}
\end{align}
Consider the integral operator
\begin{equation}\label{3.28}
(\mathcal K^{(\alpha)}f)(x)=\rho(x)^\alpha
v(x)^{1-\alpha}\int_x^\infty u(t)^{1-\alpha}f(t)dt,\quad
x\in\mathbb R.
\end{equation}
Let us show that there exists $\alpha_0\in[0,1)$ such that under
the condition $B<\infty$, the operator $\mathcal K^{(\alpha)}:
L_p\to L_p$ is bounded.

\begin{remark}\label{rem3.3} Throughout the sequel, including \S4,
we denote by $c,$ $c(p)$ some absolute positive constants which
are not essential for exposition and may differ within a single
chain of calculations.
\end{remark}
Below, when estimating $\|\mathcal K^{(\alpha)}\|_{p\to p},$ we
consecutively use Theorems \ref{thm2.13} and \ref{thm2.1}:
\begin{align}
&\|\mathcal K^{(\alpha)}\|_{p\to p}\le c(p)\sup_{x\in\mathbb
R}\left[\int_{-\infty}^x( \rho(t)^\alpha
v(t)^{1-\alpha})^pdt\right]^{1/p}\left[\int_x^\infty(u(t)^{1-\alpha})^{p'}dt\right]
^{1/p'}\nonumber\\
&\quad=c(p)\sup_{x\in\mathbb
R}\left[\int_{-\infty}^x\rho(t)^{\frac{1+\alpha}{2}p}\exp\left(-\frac{1-\alpha}{2}p
\int_t^x\frac{d\xi}{r(\xi)\rho(\xi)}\right)\exp\left(\frac{1-\alpha}{2}p\int_{x_0}^x
\frac{d\xi}{r(\xi)\rho(\xi)}\right)dt\right]^{1/p}\nonumber\\
&\quad\quad \cdot
\left[\int_x^\infty\rho(t)^{\frac{1-\alpha}{2}p'}\exp\left(-\frac{1-\alpha}{2}p'\int
_x^t\frac{d\xi}{r(\xi)\rho(\xi)}\right)\exp\left(-\frac{1-\alpha}{2}p'\int_{x_0}^x\frac
{d\xi}{r(\xi)\rho(\xi)}\right)dt\right]^{1/p'}\nonumber\\
 &\quad =c(p)\sup_{x\in\mathbb R}
 \left[\int_{-\infty}^x\rho(t)^{\frac{1+\alpha}{2}p}\exp\left(-\frac{1-\alpha}{2}p
 \int_t^x\frac{d\xi}{r(\xi)\rho(\xi)}\right)dt\right]^{1/p}\nonumber\\
 &\quad\quad
 \cdot\left[\int_x^\infty\rho(t)^{\frac{1-\alpha}{2}p'}\exp\left(-\frac{1-\alpha}{2}p'
 \int_x^t\frac{d\xi}{r(\xi)\rho(\xi)}\right)dt\right]^{1/p'}.\label{3.29}
\end{align}
Let $\alpha_0$ be a solution of the equation
$$\frac{1+\alpha}{2}p=\frac{1-\alpha}{2}p'\quad\Rightarrow\quad
\alpha_0=\frac{p'-p}{p'+p}\in[0,1).$$ Then for $\alpha=\alpha_0$
the estimate \eqref{3.29} can be simplified to the following form:
\begin{align}
\|\mathcal K^{(\alpha_0)}\|_{p\to p}&\le c(p)\sup_{x\in\mathbb
R}\left[\int_{-\infty}^x\rho(t)\exp\left(-(p-1)\int_t^x\frac{d\xi}
{r(\xi)\rho(\xi)}
\right)dt\right]^{1/p}\nonumber\\
 &\quad \cdot\left[
 \int_x^\infty\rho(t)\exp\left(-\int_x^t\frac{d\xi}{r(\xi)\rho(\xi)}\right)dt
 \right]^{1/p'}.\label{3.30}
\end{align}
Below we continue estimate \eqref{3.30} and consecutively apply
Lemmas \ref{lem2.10}, \ref{lem2.8} and \ref{lem2.5}:
\begin{align}
\|\mathcal K^{(\alpha_0)}\|_{p\to p}&\le c(p)\sup_{x\in\mathbb
R}\left[\sum_{n=-\infty}^{-1}\int_{\Delta_n}\rho(t)\exp\left(-(p-1)\int_t^{\Delta_{-1}^+}
\frac{d\xi}{r(\xi)\rho(\xi)}\right)dt\right]^{1/p}\nonumber\\
&\quad
\cdot\left[\sum_{n=1}^\infty\int_{\Delta_n}\rho(t)\exp\left(-\int_{\Delta_1^-}^t
\frac{d\xi}{r(\xi)\rho(\xi)}\right)dt\right]^{1/p'}\nonumber\\
&\le c(p)\sup_{x\in\mathbb
R}\left[\sum_{n=-\infty}^{-1}\rho(x_n)d(x_n)\exp\left(-(p-1)\int_{\Delta_n^+}^
{\Delta_{-1}^+}\frac{d\xi}{r(\xi)\rho(\xi)}\right)\right]^{1/p}\nonumber\\
&\quad\cdot
\left[\sum_{n=1}^\infty\rho(x_n)d(x_n)\exp\left(-\int_{\Delta_1^-}^{\Delta_n^-}
\frac{d\xi}{r(\xi)\rho(\xi)}\right)\right]^{1/p'}\nonumber\\
&\le c(p)\sup_{x\in\mathbb
R}\left[\sum_{n=-\infty}^{-1}h(x_n)d(x_n)\exp\left(-\frac{p-1}{2}\int_{\Delta_n^+}^{
\Delta_{-1}^+} \frac{d\xi}{r(\xi)h(\xi)}\right)\right]^{1/p}\nonumber\\
&\quad \cdot\left[\sum_{n=1}^\infty
h(x_n)d(x_n)\exp\left(-\frac{1}{2}\int_{\Delta_1^-}^{\Delta_n^-}\frac{d\xi}{
r(\xi)h(\xi)}\right)\right]^{1/p'}\label{3.31}.
\end{align}

 Note that Definition \ref{defn2.9} and \eqref{1.18} imply the
 equalities
  \begin{equation}\label{3.32}
\begin{cases}
\displaystyle{\int_{\Delta_n^+}^{\Delta_{-1}^+}\frac{d\xi}{r(\xi)h(\xi)}}
=|n|-1\quad
 &\text{if}\quad n\le-1\\ \\
 \displaystyle{\int_{\Delta_1^-}^{\Delta_{n}^-}\frac{d\xi}{r(\xi)h(\xi)}} =n-1\quad
 &\text{if}\quad n\ge1
 \end{cases}
 \end{equation}
 Equalities \eqref{3.32} are checked in the same way, and
 therefore we only check, say, the second one. For $n=1$ it is
 obvious, and for $n\ge2$ we have
 $$\int_{\Delta_1^-}^{\Delta_n^-}\frac{d\xi}{r(\xi)h(\xi)}=\sum_{k=1}^{n-1}\int
 _{\Delta_k}\frac{d\xi}{r(\xi)h(\xi)}=\sum_{k=1}^{n-1}1=n-1\quad\Rightarrow\quad
 \eqref{3.32}.$$

 Let us now continue estimate \eqref{3.31} using \eqref{3.32}:
 \begin{align*}
 \|\mathcal K^{(\alpha_0)}\|_{p\to p}&\le c(p)\sup_{x\in\mathbb
 R}\left[\sum_{n=-\infty}^{-1}h(x_n)d(x_n)\exp\left(-\frac{p-1}{2}(|n|-1)\right)
 \right]^{1/p} \\
 &\quad\cdot\left[\sum_{n=1}^\infty
 h(x_n)d(x_n)\exp\left(-\frac{n-1}{2}\right)\right]^{1/p'} \\
 &\le
 c(p)B^{\frac{1}{p}+\frac{1}{p'}}\left[\sum_{n=-\infty}^{-1}\exp\left(-\frac{
 (p-1)(|n|-1)}{2}\right)\right]^{1/p}\\
 &\quad\cdot\left[\sum_{n=1}^\infty\exp\left(-
 \frac{n-1}{2}\right)\right]^{1/p'}=c(p)B\quad\Rightarrow
 \end{align*}
 \begin{equation}\label{3.33}
\|\mathcal K^{(\alpha_0)}\|_{p\to p}\le c(p)B\quad\text{for}\quad
 p\in(1,2],\qquad \alpha_0=\frac{p'-p}{p'+p}.
\end{equation}

Let now $p\in(2,\infty),$\ $f\in L_p,$ and denote by $\alpha $
some number from the interval $(0,1)$ which will be chosen later.
Then by \lemref{lem1.3}, for any $x\in\mathbb R,$ we get
\begin{align}
|(G_2f)(x)|&\le v(x)\int_x^\infty u(t)|f(t)|dt\le
v(x)^{1-\alpha}\int_x^\infty(v(t)u(t))^\alpha
u(t)^{1-\alpha}|f(t)|dt\nonumber\\
&=v(x)^{1-\alpha}\int_x^\infty\rho(t)^\alpha
u(t)^{1-\alpha}|f(t)|dt.\label{3.34}
\end{align}
Furthermore, in the same way as above, for the norm of the
operator
$$(\tilde{\mathcal
K}^{(\alpha)}f)(x)=v(x)^{1-\alpha}\int_x^\infty\rho(t)^\alpha
u(t)^{1-\alpha}f(t)dt,\quad x\in\mathbb R$$ for
$\alpha=\alpha_1=\frac{p-p'}{p+p},$ we establish the estimate
\begin{equation}\label{3.35}
 \|\tilde{\mathcal K}^{(\alpha_0)}\|_{p\to p}\le c(p)B,\quad
p\in(2,\infty).
\end{equation}
From \eqref{3.27} and \eqref{3.33},\ \eqref{3.34} and
\eqref{3.35}, it now follows that
$$\|G_2f\|_p\le
\| {\mathcal K}^{(\alpha_0)}|f|\,\|_p\le\| {\mathcal
K}^{(\alpha_0)}\|_{p\to p}\|f\|_p\le c(p)B\|f\|_p,\quad
p\in(1,2];$$
$$\|G_2f\|_p\le
\| \tilde{\mathcal K}^{(\alpha_1)}|f|\,\|_p\le\|\tilde {\mathcal
K}^{(\alpha_1)}\|_{p\to p}\|f\|_p\le c(p)B\|f\|_p,\quad
p\in(2,\infty).$$ These estimates imply inequalities \eqref{3.26}.
This completes the proof of the theorem.
\end{proof}

\begin{proof}[Proof of Corollary \ref{cor1.9}] From conditions
\eqref{1.2} and \eqref{1.6} it follows that all auxiliary
functions from \lemref{lem1.6} are well-defined and, in addition,
condition \eqref{1.5} holds. From the properties of the FSS
$\{u,v\}$ of equation \eqref{1.7} (see \lemref{lem1.3}) for all
$x\in\mathbb R$, we obtain the inequalities
\begin{equation}\label{3.36}
r(x)v'(x)\ge\int_{-\infty}^x q(t)v(t)dt,\qquad
r(x)|u'(x)|\ge\int_x^\infty q(t)u(t)dt.
\end{equation}

Below we consecutively use equation \eqref{1.9}, estimates
\eqref{3.36}, formula \eqref{1.12}, \lemref{lem1.6},
representation \eqref{2.3}, inequalities \eqref{2.9} and
\eqref{2.6}, equality \eqref{1.18} and the definition \eqref{1.20}
of $A:$
\begin{align*}
1&=r(x)v'(x)u(x)-r(x)u'(x)v(x)\ge u(x)\int_{-\infty}^x
q(t)v(t)dt+v(x)\int_x^\infty q(t)u(t)dt\\
&=\int_{-\infty}^\infty
q(t)G(x,t)dt\ge\int_{x-d(x)}^{x+d(x)}q(t)G(x,t)dt\\
&=\int_{x-d(x)}^{x+d(x)}q(t)\sqrt{\rho(t)\rho(x)}\exp\left(-\frac{1}{2}\left|
\int_x^t\frac{d\xi}{r(\xi)\rho(\xi)}\right|\right)dt\\
&\ge
c^{-1}\rho(x)\exp\left(-\frac{1}{2}\int_{x-d(x)}^{x+d(x)}\frac{d\xi}{r(\xi)\rho(\xi)}\right)
\int_{x-d(x)}^{x+d(x)}q(t)dt\\
&\ge
c^{-1}\exp\left(-\int_{x-d(x)}^{x+d(x)}\frac{d\xi}{r(\xi)h(\xi)}\right)h(x)d(x)
\left[\frac{1}{2d(x)}\int_{x-d(x)}^{x+d(x)}q(t)dt\right]\ge
c^{-1}A(h(x)d(x))\quad\Rightarrow\\
&\qquad\qquad h(x)d(x)\le cA^{-1}<\infty,\quad x\in\mathbb
R\quad\Rightarrow\quad B=\sup_{x\in\mathbb R}h(x)d(x)\le
cA^{-1}<\infty.
\end{align*}
The assertion of the Corollary now follows from \thmref{thm1.8}.
\end{proof}

\renewcommand{\qedsymbol}{}
\begin{proof}[Proof of Corollary \ref{cor1.10}]  Necessity.

Suppose equation \eqref{1.1} is correctly solvable in $L_p,$ \
$p\in(1,\infty).$ Since $r\equiv1,$ equalities \eqref{1.4} hold,
and by \thmref{thm1.1} relations \eqref{1.5} and \eqref{1.6} are
satisfied. Then all auxiliary functions from \lemref{lem1.6} and
the function $\tilde d(x)$ from \lemref{lem2.6} are well-defined
(see \cite{7}). Furthermore, from \eqref{1.18}, \eqref{2.6} and
\eqref{2.9}, we get
\begin{align*}
\begin{cases} 1=\displaystyle{\int_{x-d(x)}^{x+d(x)}\frac{dt}{h(t)}\ge
\frac{1}{2}\int_{x-d(x)}^{x+d(x)}\frac{dt}{\rho(t)}
\ge\frac{1}{2e^2}\ \frac{d(x)}{\rho(x)}},&\quad x\in\mathbb R\\
&\quad\qquad\qquad \quad \Rightarrow\\
 1=\displaystyle{\int_{x-d(x)}^{x+d(x)}\frac{dt}{h(t)}\le
2\int_{x-d(x)}^{x+d(x)}\frac{dt}{\rho(t)} \le {2e^2}\
\frac{d(x)}{\rho(x)}},&\quad x\in\mathbb R
\end{cases}
\end{align*}

\begin{equation}\label{3.37}
2^{-1}e^{-2}d(x)\le \rho(x)\le 2e^2d(x),\quad x\in\mathbb R.
\end{equation}

{}From \eqref{3.37}, \eqref{2.6} and \eqref{2.8}, it is easy to
obtain the estimates
\begin{align}
 c^{-1}\tilde d(x)\le h(x),d(x)\le c \tilde d(x),&\quad  x\in\mathbb
R\quad\Rightarrow\nonumber\\
 c^{-1}\tilde d(x)^2\le h(x) d(x)\le c \tilde d(x)^2,&\quad
x\in\mathbb R\label{3.38}
\end{align}
Since $B<\infty$ in view of \thmref{thm1.8}, according to
\eqref{3.38}, we have $d_0=\sup\limits_{x\in\mathbb R}\tilde
d(x)<\infty.$ We now obtain from \eqref{2.7}:
$$2=\tilde d(x)\int_{x-\tilde d(x)}^{x+\tilde d(x)}q(t)dt\le
d_0\int_{x-d_0}^{x+d_0}q(t)dt\quad\Rightarrow\quad\inf_{x\in\mathbb
R}\int_{x-d_0}^{x+d_0}q(t)dt\ge\frac{2}{d_0}>0.$$
\end{proof}

\renewcommand{\qedsymbol}{\openbox}
 \begin{proof}[Proof of \thmref{thm1.4}]  Sufficiency.

Suppose \eqref{1.21} holds. Then it is clear that relations
\eqref{1.5} and \eqref{1.6} are satisfied, and therefore all
auxiliary functions from \lemref{lem1.6} and the function $\tilde
d$ from \lemref{lem2.6} are well-defined.  Hence inequalities
\eqref{3.38} remain true. From \eqref{1.21} it follows that there
exists $d_0\gg1$ such that
$$\int_{x-d_0}^{x+d_0}q(t)dt>\frac{2}{d_0},\quad x\in\mathbb R.$$
Then $\tilde d(x)\le d_0$ for every $x\in\mathbb R.$ Indeed, if
for some $x\in\mathbb R$ this inequality does not hold, then we
have a contradiction:
$$2=\tilde d(x)\int_{x-\tilde
d(x)}^{x+\tilde d(x)}q(t)dt\ge d_0\int_{x-d_0}^{x+d_0} q(t)dt>2.$$
Thus $\tilde d(x)\le d_0$ for every $x\in\mathbb R.$ But then from
\eqref{3.38} we obtain $B\le cd_0^2<\infty.$ It remains to refer
to \thmref{thm1.8}.
\end{proof}

\section{Additional assertions and examples}

Below we present several   applications of the results obtained
above.

\begin{thm}\label{thm4.1} \cite{2} Suppose conditions \eqref{1.2}
hold, and, in addition, $q_0>0,$ where
\begin{equation}\label{4.1}
q_0=\inf_{x\in\mathbb R} q(x).
\end{equation}
Then equation\eqref{1.1} is correctly solvable in $L_p$ for all
$p\in(1,\infty).$
\end{thm}

\begin{proof}
Since $q_0>0$, condition \eqref{1.6} holds, and by \lemref{lem1.6}
the function $d$ is well-defined. Then (see \eqref{1.20}):
$$A=\inf_{x\in\mathbb
R}\frac{1}{2d(x)}\int_{x-d(x)}^{x+d(x)}q(t)dt\ge\inf_{x\in\mathbb
R}\frac{1}{2d(x)}\int_{x-d(x)}^{x+d(x)}q_0dt=q_0>0.$$ The
statement of the theorem now follows from Corollary \ref{cor1.9}
and Remark \ref{rem1.11}.
\end{proof}

\begin{remark}\label{rem4.2} The proof of \thmref{thm4.1} given
here was proposed in \cite{7}.
\end{remark}

Thus, the requirement $q_0>0$ to the function $q$ is so strong
that the answer to the question on I)--II) does not depend on the
behavior of the function $r$ (assuming \eqref{1.2}). In this
connection, consider the opposite direction, i.e., find
requirements to the function $r$ under which the solution of the
problem I)--II) does not depend on the behavior of the function
$q$ (in a certain framework).

\begin{thm}\label{thm4.3} Suppose conditions \eqref{1.2} hold,
and, in addition,
\begin{equation}\label{4.2}
\int_{-\infty}^0q(t)dt=\int_0^\infty q(t)dt=\infty.
\end{equation}
Then equation \eqref{1.1} is correctly solvable in $L_p$ for all
$p\in(1,\infty)$ if $\theta<\infty.$ Here
\begin{equation}\label{4.3}
\theta=\sup_{x\in\mathbb
R}|x|\left(\int_{-\infty}^x\frac{dt}{r(t)}\right)\left(\int_x^\infty\frac{dt}{r(t)}\right).
\end{equation}
\end{thm}

\begin{proof}
{}From \eqref{1.2} and \eqref{4.2} it follows that relations
\eqref{1.5} and \eqref{1.6} hold, and this means that the
assumptions of Lemmas \ref{lem1.3} and \ref{lem1.6} and
\thmref{thm1.8} are satisfied. Furthermore, from \lemref{lem1.3}
it is easy to obtain the relations
\begin{gather}
u(x)=v(x)\int_x^\infty\frac{dt}{r(t)v^2(t)},\quad
v(x)=u(x)\int_{-\infty}^x\frac{dt}{r(t)u^2(t)},\quad x\in\mathbb
R\quad\Rightarrow\nonumber\\
\rho(x)=u(x)v(x)=v^2(x)\int_x^\infty\frac{dt}{r(t)v^2(t)}=u^2(x)\int_{-\infty}^x
\frac{dt}{r(t)u^2(t)},\quad x\in\mathbb R\label{4.4}
\end{gather}
Note that $\frac{1}{r}\in L_1$ because $\theta<\infty.$ Therefore
from \eqref{4.4} and \lemref{lem1.3}, it is easy to obtain the
following conceptual estimates:
\begin{equation}
\begin{aligned}\label{4.5}
\rho(x)\le
u^2(x)\int_{-\infty}^x\frac{dt}{r(t)u^2(x)}=\int_{-\infty}^x
\frac{dt}{r(t)},&\quad\text{if}\quad x\le 0\\
\rho(x)\le v^2(x)\int_{x
}^\infty\frac{dt}{r(t)v^2(x)}=\int_{x}^\infty
\frac{dt}{r(t)},&\quad\text{if}\quad x\ge 0.
\end{aligned}
\end{equation}

We shall also need the following simple consequence of
\lemref{lem1.6}:
\begin{equation}\label{4.6}
d(x)\le|x|\qquad\text{for all}\qquad |x|\gg 1.
\end{equation}

Below we consecutively use \lemref{lem2.5}, \eqref{4.6},
\eqref{4.5} and the condition
\begin{equation}\label{4.7}
h(x)d(x)\le 2\rho(x)x\le 2x\int_x^\infty\frac{dt}{r(t)}\le
c<\infty\quad\text{for}\quad x\gg1;
\end{equation}
\begin{equation}\label{4.8}
h(x)d(x)\le 2\rho(x)|x|\le 2|x|\int_{-\infty}^x\frac{dt}{r(t)}\le
c<\infty\quad\text{for}\quad x\ll-1.
\end{equation}
Moreover, since the function $\rho(x)d(x)$ is continuous and
positive for $x\in\mathbb R$ (see \lemref{lem1.6}), the
inequalities
$$0< h(x)d(x)\le 2\rho(x)d(x),\qquad x\in\mathbb R$$
imply that the function $h(x)d(x)$ is bounded on every finite
interval. Together with \eqref{4.7} and \eqref{4.8}, this give
$B<\infty$ (see \eqref{1.19}). It remains to apply
\thmref{thm1.8}.
\end{proof}

The following applications of \thmref{thm4.3} to concrete
equations are subdivided into pairs so that it is interesting to
compare the examples constituting each pair.

\begin{example}\label{examp4.4} Consider equation \eqref{1.1} with
coefficients
$$ r(x)=\begin{cases}1&\ \text{if}\quad |x|\le 1\\
&\qquad\qquad\qquad ,\ q(x)=1+\cos|x|^\alpha,\quad x\in\mathbb
R,\quad \alpha\in(0,\infty)\\
x^2&\ \text{if}\quad |x|\ge 1
\end{cases}$$
In this case the hypothesis of \thmref{thm4.3} are obviously
satisfied, $\theta<\infty$ and therefore such an equation is
correctly solvable in $L_p$ for all $p\in(1,\infty)$ regardless of
$\alpha\in(0,\infty).$
\end{example}

\begin{example}\label{examp4.5} Consider equation \eqref{1.1} with
coefficients
$$ r(x)\equiv1,\quad x\in\mathbb R;\qquad q(x)=1+\cos|x|^\alpha,\qquad
\alpha\in(0,\infty),\quad x\in\mathbb R.$$ Then (see \cite{6})
such an equation is correctly solvable in $L_p$  if and only if
$\alpha\ge1$.
\end{example}

\begin{remark}\label{rem4.6}
The statement of Example \ref{examp4.5} becomes completely obvious
if one compares the criterion for correct solvability $m(a)>0,$ \
$a\in(0,\infty)$ (see \eqref{1.21}) with the different behavior at
infinity of the graphs of $q$ for $\alpha\in(0,1)$ and
$\alpha\in[1,\infty)$ in the zeros of the function $q.$
\end{remark}

Let us now consider two examples of a different type.

\begin{example}\label{examp4.7} Consider equation \eqref{1.1} with
coefficients
\begin{equation}\label{4.9}
 r(x)=\begin{cases} 1,&\ \text{if}\quad  |x|\le 1\\
x^2,&\ \text{if}\quad |x|\ge 1
\end{cases},\quad q(x)=\begin{cases} 1,&\ \text{if}\quad |x|\le 1\\ \frac{1}{\sqrt
{|x|}},& \ \text{if}\quad |x|\ge 1
\end{cases}
\end{equation}
In the case \eqref{4.9} all the hypotheses of \thmref{thm4.3} are
obviously satisfied, $\theta<\infty$ and therefore such an
equation is correctly unsolvable in $L_p$ for all
$p\in(1,\infty).$
\end{example}

\begin{example}\label{examp4.8} Consider equation \eqref{1.1} with
coefficients $r\equiv1$,\ $x\in\mathbb R$ and $q$ satisfying
condition \eqref{1.2}. If, in addition,
\begin{equation}\label{4.10}
\lim_{x\to-\infty}q(x)=0\qquad\text{or}\qquad
\lim_{x\to\infty}q(x)=0,
\end{equation}
equation \eqref{1.1} is correctly unsolvable in $L_p$ for all
$p\in(1,\infty)$ (see Corollary \ref{cor1.10}).

\end{example}

Let us now consider the equation of direct application of
\thmref{thm1.8} to particular  equations \eqref{1.1}. Since it is
usually impossible to find an explicit form for the functions $h$
and $d$, a general method of applying the theorem consists of
obtaining two-sided, sharp by order estimates for the functions
$h$ and $d$ at infinity. Clearly, such inequalities lead to exact
estimates of the function $h\cdot d$ at infinity and then to a
complete answer to the question on the finiteness of $B$ (see the
proof of \thmref{thm4.3}) and, finally, to a concluding statement
using \thmref{thm1.8}.

Below we present an example where this scheme is realized for the
case of equation \eqref{1.1} with coefficients \eqref{4.9}. To
this end, we present the following assertion.

\begin{lem}\label{lem4.9} Suppose conditions \eqref{1.2} and
\eqref{1.6} hold. For a given $x\in\mathbb R$ let us introduce the
functions $F_1(\eta),$ $F_2(\eta)$ and $F_3(\eta)$ with
$\eta\ge0:$
\begin{gather*}
F_1(\eta)=\int_{x-\eta}^x\frac{dt}{r(t)}\cdot \int_{x-\eta}^x
q(t)dt,\qquad F_2(\eta)=\int_x^{x+\eta}\frac{dt}{r(t)}\cdot
\int_x^{x+\eta}q(t)dt,\\
F_3(\eta)=\int_{x-\eta}^{x+\eta}\frac{dt}{r(t)h(t)}.
\end{gather*}
Then the following assertions hold (see \lemref{lem1.6}):
\begin{enumerate}
\item[a)] the inequality $\eta\ge d_1(x)\ (0\le\eta\le d_1(x))$
holds if and only if $F_1(\eta)\ge1$\ $(F_1(\eta)\le 1);$
\item[b)] the inequality $\eta\ge d_2(x)\ (0\le\eta\le d_2(x))$
holds if and only if $F_2(\eta)\ge1$\ $(F_2(\eta)\le 1);$
\item[c)] the inequality $\eta\ge d(x)\ (0\le\eta\le d (x))$ holds
if and only if $F_3(\eta)\ge1$\ $(F_3(\eta)\le 1).$
\end{enumerate}
\end{lem}

\renewcommand{\qedsymbol}{}
\begin{proof}[Proof of \lemref{lem4.9}]  Necessity.

All assertions of the Lemma are proved in the same way. Consider,
say, b). Clearly,
\begin{equation}\label{4.11}
F_2'(\eta)=\frac{1}{r(x+\eta)}\int_x^{x+\eta}
q(t)dt+q(x+\eta)\int_x^{x+\eta}\frac{dt}{r(t)},
\end{equation}
and therefore $F_2'(\eta)\ge0$ in view of \eqref{1.2}. Then if
$\eta\ge d_2(x),$ then $F_2(\eta)\ge F_2(d_2(x))=1.$
\end{proof}

\renewcommand{\qedsymbol}{\openbox}
 \begin{proof}[Proof of \lemref{lem4.9}]  Sufficiency.
Let $F_2(\eta)\ge1.$ Assume the contrary: $\eta<d_2(x).$ Since
$F_2(\eta)\ge1,$ we have
$$\int_x^{x+\eta}q(t)dt>0\quad\Rightarrow\quad F_2'(\eta)>0\quad
\text(see\ \eqref{4.11})\quad\Rightarrow\quad 1\le
F_2(\eta)<F_2(d_2(x))=1,$$ contradiction. Hence $\eta\ge d_2(x).$
 \end{proof}

 \begin{remark}\label{rem4.10}
 \lemref{lem4.9} is an efficient tool for proving estimates of the
 functions $h$ and $d$ (see below). In particular, in \cite{10} it
 was used for a meaningful class of equations \eqref{1.1} in order
 to get a priori, sharp by order, two-sided estimates of these
 functions expressed in terms of the functions $r$ and $q.$ See
 \cite{7,10} for a detailed exposition of proofs and applications
 of such inequalities. In the case \eqref{4.9}, a priori estimates
 from \cite{10} are not applicable because of the fast growth of
 the function $r/q.$
\end{remark}

 Let us introduce the following notation. Let $\alpha(x)$ and
 $\beta(x)$ be continuous positive functions for $x\in (a,b) $\ $(-\infty\le a
 <b\le\infty)$.
 We write $\alpha(x)\asymp\beta(x)$ for
 $x\in(a,b)$ if there exists a constant $c\in[1,\infty)$ such that
 $$c^{-1}\alpha(x)\le\beta(x)\le c\alpha(x)\qquad\text{for
 all}\qquad x\in(a,b).$$

 \begin{lem}\label{lem4.11}
 In the case \eqref{4.9}, the following relations hold (see
 \lemref{lem1.6}):
 \begin{equation}\label{4.12}
 \varphi(x)\asymp\begin{cases}\frac{1}{|x|}\quad &\text{if}\quad
 x\ll-1\\
 \frac{1}{\sqrt{|x|}}\quad &\text{if}\quad
 x\gg 1\end{cases}\qquad \psi(x)\asymp\begin{cases}\frac{1}{\sqrt{|x|}}
 \quad &\text{if}\quad
 x\ll-1\\
 \frac{1}{ |x|}\quad &\text{if}\quad
 x\gg 1\end{cases}
 \end{equation}
 \begin{equation}\label{4.13}
 h(x)\asymp\frac{1}{|x|}\qquad\text{if}\qquad |x|\gg 1.
 \end{equation}
\end{lem}

\begin{proof}
Both relations of \eqref{4.12} are proved in the same way. Let us
check, say, the second one. It is easy to see that in the case
\eqref{4.9} the following equality holds (see \lemref{lem4.9}):
\begin{equation}\label{4.14}
F_2(\eta)=\frac{2\eta^2}{x(x+\eta)(\sqrt{|x|}+\sqrt{|x+\eta|)}},\quad{if}\quad
[x,x+\eta]\cap[-1,1]=\emptyset.
\end{equation}
Let $x\gg1.$ Then in view of \eqref{4.14} we have the inequalities
$$F_2(\eta)\bigm|_{\eta=4x^2}>1,\qquad
F_2(\eta)\bigm|_{\eta=\frac{x^2}{4}}<1,$$ and therefore by
\lemref{lem4.9}, we obtain
$$4^{-1}x^2\le d_2(x)\le 4x^2\quad\text{for}\quad x\gg 1.$$

Similarly, for $x\ll-1$ in view of \eqref{4.14}, we have the
inequalities: $$F_2(\eta)\bigm|_{\eta=|x|-\sqrt{|x|}}>1,\quad
 F_2(\eta)\bigm|_{\eta=|x|-4 \sqrt{|x|}}<1,
$$ and therefore by
\lemref{lem4.9} we obtain:
$$|x|-4\sqrt{|x|}\le d_2(x)\le
|x|-\sqrt{|x|}\qquad\text{for}\qquad x\ll-1.$$ Thus,
\begin{equation}\label{4.15}
d_2(x)\asymp\begin{cases} x^2&\ \text{if}\quad x\gg1\\
|x|&\ \text{if}\quad x\ll-1\end{cases}
\end{equation}

Relations \eqref{4.12} for $\psi$ follows from \eqref{4.15}, the
definition of the function $\psi$ and \eqref{4.9} Formula
\eqref{4.13} follows from \eqref{4.12} and the definition of the
function $h.$ \end{proof}

\begin{remark}\label{rem4.12} If one follows the main method  of \lemref{lem4.11}, then using
\eqref{4.13} and \lemref{lem4.9}  one can obtain two-sided, sharp
by order estimates for the function $d.$ However, it is worth
noting that in this situation, as in many others, one can always
``economize" technical work if the already obtained estimate for
the function $h$ shows that $h(x)|x|\le c<\infty$ for all
$|x|\gg1.$ Indeed, in the latter case the proof of the estimates
for the function $d$ becomes superfluous because one can replace
the sharp by order upper estimate of the function $d$ with the
rougher a priori estimate \eqref{4.6} without changing the results
on $B$ (see \thmref{thm1.8}).
\end{remark}

Let us now show, say, that $B<\infty$ (see \eqref{1.19}). From
\eqref{4.13} and \eqref{4.6}, we get the inequalities
\begin{equation}\label{4.16}h(x)d(x)\le\frac{c}{|x|}\cdot |x|=c<\infty\qquad\text{for
all}\qquad |x|\gg1.\end{equation}
 Since on every finite segment
$[a,b]$ the function $h(x)d(x),$ \ $x\in[a,b]$ is bounded, from
\eqref{4.16} we conclude that indeed $B<\infty.$ By
\thmref{thm1.8}, this implies that in the case \eqref{4.9},
equation \eqref{1.1} is correctly solvable in $L_p$ for all
$p\in(1,\infty).$


\begin{thebibliography}{10}

\bibitem[1]{3} R.C. Brown and J. Cook, \emph{Continuous invertibility of minimal
Sturm-Liouville operators in Lebesgue spaces},
  Proc. of the Royal Soc. of Edinburgh \textbf{136A}  (2006), 53-70.

 \bibitem[2]{13} N. Chernyavskaya and L. Shuster, \emph{Weight summability of solutions of the
 Sturm-Liouville equation},
  J. Diff. Eq. \textbf{151}  (1999), 456-473.

\bibitem[3]{7} N. Chernyavskaya and L. Shuster, \emph{Estimates for the Green function
 of a general Sturm-Liouville
operator and their  applications}, Proc. Amer. Math. Soc.
\textbf{127} (1999), 1413-1426.


 \bibitem[4]{10} N. Chernyavskaya and L. Shuster, \emph{Regularity of the inversion
 problem for a Sturm-Liouville equation in $L_p(R)$},  Methods and Applications of
 Analysis \textbf{7}
  (2000), 65-84.

 \bibitem[5]{6} N. Chernyavskaya and L. Shuster, \emph{A criterion for correct
 solvability of the  Sturm-Liouville equation in the space $L_p(R)$}, Proc. Amer. Math.
 Soc. \textbf{130}
  (2001), no.~4,  1043-1054.

\bibitem[6]{17} N. Chernyavskaya and L. Shuster, \emph{Conditons for correct solvability
of a simplest singular boundary value problem of general form.
I.}, Z. Anal. Anwendungen \textbf{25}
  (2006), 205-235.

  \bibitem[7]{1} R. Courant \emph{Partial Differential Equations},
 John Wiley \& Sons,
  New York, 1962.


 \bibitem[8]{16} E.B. Davies and E.M. Harrell, \emph{Conformally flat Riemannian metrics, Schr\"odinger
 operators and semiclassical approximation},
J. Diff. Eq. \textbf{66} (1987), 165-188.

\bibitem[9]{14} W.N. Everitt and M. Giertz, \emph{Some properties of the domains of
certain differential operators}, Proc. London Math. Soc.
\textbf{23} (1971) (3), 301-324.

\bibitem[10]{12} W.N. Everitt and M. Giertz, \emph{Some inequalities associated with
certain differential operators}, Math. Z.. \textbf{126} (1972)
(4), 308-326.

\bibitem[11]{19} L.W. Kantorovich and G.P. Akilov, \emph{Functional Analysis}, Nauka, Moscow, 1977.

\bibitem[12]{9} K. Mynbaev and M. Otelbaev, \emph{Weighted Function Spaces and the
Spectrum of Differential Operators}, Nauka, Moscow, 1988.

\bibitem[13]{2} R. Oinarov, \emph{Properties of Sturm-Liouville operator in $L_p$},
Izvestija Akad. Nauk Kazakh. SSR \textbf{1} (1990), 43-47.

\bibitem[14]{8} R. Oinarov and M. Otelbaev, \emph{A criterion for discreteness of the
spectrum of the general Sturm-Liouville operator, and embedding
theorems connected with it}, Diff. Uravn. \textbf{24} (1988),
584-591; (Engl. transl. Diff. Eqns. \textbf{24} (1988), 402-408).


\bibitem[15]{11} B. Opic and A. Kufner, \emph{Hardy Type Inequalities},
Pitman Research Notes in Mathematics Series 219, Harlow, Longman,
1990.


\bibitem[16]{15} M. Otelbaev, \emph{A criterion for the resolvent of a Sturm-Liouville operator to
be a kernel}, Math. Notes \textbf{25} (1979), 296-297 (translation
of Mat. Zametki).

\bibitem[17]{18} M. Otelbaev, \emph{Estimates of Spectrum of Sturm-Liouville Operator},
 (in Russian), Alma-Ata, Gilim, 1990.

 \bibitem[18]{4} M. Sapenov and L. Shuster, \emph{Estimates of Green's function and theorem
 on separability of a Sturm-Liouville operator in $L_p$},
Manuscript deposited at VINITI (1985), No.~8257-B85.

\bibitem[19]{5} L. Shuster, \emph{A priori properties of solutions of a Sturm-Liouville
equation and A.M. Molchanov's criterion}, Math. Notes Ac. Sci.
USSR \textbf{50} (1991), 746-751; (translation of Mat. Zametki).

  \end{thebibliography}
\end{document}